\overfullrule=0pt
\centerline {\bf Kirchhoff-type problems involving subcritical and superlinear nonlinearities satisfying no further condition}\par
\bigskip
\bigskip
\centerline {BIAGIO RICCERI}\par
\bigskip
\bigskip
\centerline {\it Dedicated to the memory of Anna Aloe}\par
\bigskip
\bigskip
{\bf Abstract.} In this note, we deal with a problem of the type
$$\cases {-h\left ( \int_{\Omega}|\nabla u(x)|^2dx\right )
\Delta u=f(u) & in $\Omega$\cr & \cr
u_{|\partial\Omega}=0\ .\cr}$$
As an application of a new general multiplicity result, we establish the existence of at least three solutions, two of which are global minima of the related energy functional. The only condition assumed on $f$ is that it be subcritical and superlinear: no condition on the behaviour of $f$ at $0$ is requested.
\bigskip
Math Subject Classifications: 35J20, 35J61, 49K40, 90C26.\par
\bigskip
Key words: Kirchhoff-type problems; multiplicity of global minima; variational methods.\par
\bigskip
\bigskip
\bigskip
\bigskip
Here and in what follows, $\Omega\subset {\bf R}^m$ is a smooth bounded domain, with $m\geq 3$.
\smallskip
For every $q\in \left ] 0, {{m+2}\over {m-2}}\right ]$, we denote by ${\cal A}_q$ the class of all continuous 
functions $f: {\bf R}\to {\bf R}$ such that
$$\limsup_{|\xi|\to +\infty}{{|f(\xi)|}\over {|\xi|^q}}<+\infty$$
and
$$-\infty<\liminf_{|\xi|\to +\infty}{{F(\xi)}\over {\xi^2}}\leq\limsup_{|\xi|\to +\infty}{{F(\xi)}\over {\xi^2}}=+\infty$$
where $F(\xi)=\int_0^{\xi}f(t)dt\ .$\par
\smallskip
Given $f\in {\cal A}_q$ and a continuous function $h:[0,+\infty[\to {\bf R}$, 
consider the following Kirchhoff-type problem
$$\cases {-h\left ( \int_{\Omega}|\nabla u(x)|^2dx\right )
\Delta u=f(u) & in $\Omega$\cr & \cr
u_{|\partial\Omega}=0\ .\cr}$$
A weak solution of this problem is any $u\in H^1_0(\Omega)$ such that
$$h\left ( \int_{\Omega}|\nabla u(x)|^2dx\right ) \int_{\Omega}
\nabla u(x)\nabla v(x)dx=\int_{\Omega}f(u(x))v(x)dx$$
for all $v\in H^1_0(\Omega)$.\par
\medskip
So, the weak solutions of the problem are precisely the critical points in $H^1_0(\Omega)$ of the functional
$$u\to {{1}\over {2}}H\left ( \int_ {\Omega}|\nabla u(x)|^2dx\right )-\int_{\Omega}F(u(x))dx$$
where $H(t)=\int_0^th(s)ds.$\par
\smallskip
A real-valued function $g$ on a topological space is said to be sequentially inf-compact if, for each $r\in {\bf R}$, the set
$g^{-1}(]-\infty,r])$ is sequentially compact.\par
\medskip
The aim of this note is to establish the following result:\par
\medskip
THEOREM 1 . - {\it For each $q\in \left ] 0, {{m+2}\over {m-2}}\right [$ and $f\in {\cal A}_q$ 
there exists a divergent sequence $\{a_n\}$ in $]0,+\infty[$ with the following property: for
every $n\in {\bf N}$ and for every continuous and non-decreasing function
$k:[0,+\infty[\to [0,+\infty[$, with $\lim_ {t\to +\infty}{{K(t)}\over {t^{q+1\over 2}}}=+\infty$ and int$(k^{-1}(0))=\emptyset$,
there esists $b>0$ such that the problem
$$\cases{-\left (a_n+bk\left ( \int_ {\Omega}|\nabla u(x)|^2dx\right )\right )\Delta u=f(u) & in $\Omega$\cr & \cr
u_{|\partial \Omega}=0\cr}$$
has at least three weak solutions, two of which are global minima in $H^1_0(\Omega)$ of the energy functional
$$u\to {{a_n}\over {2}} \int_ {\Omega}|\nabla u(x)|^2dx+{{b}\over {2}}K\left ( \int_ {\Omega}|\nabla u(x)|^2dx\right )-
\int_{\Omega}F(u(x))dx$$
where $K(t)=\int_0^tk(s)ds\ .$}\par
\medskip
A comparison of Theorem 1 with known results cannot be properly done. This is due to the fact that no previous result on the problem we are
dealing with guarantees the existence of at least two global minima of the energy functional related to it. More precisely, no such a result is known
when the nonlinearity $f$, as in our case, does not depend on $x$ ($x\in \Omega$). For quite special $f$ depending necessarily on $x$, the only known
results of that type have been obtained in [5]. But, also for what concerns the assumptions on $f$, Theorem 1 presents a novelty: it seems that,
even when the energy functional in unbounded below, no existing result ensures the existence of at least three solutions 
of the problem assuming on $f$ only its belonging to the class ${\cal A}_q$. Actually, some condition on the behaviour of $f$ at $0$ is 
usually assumed (see, for instance, [1], [2], [4], [8], [9], [11] and references therein).
\smallskip
Our proof of Theorem 1 is based on the use of the following new abstract multiplicity result:\par
\medskip
THEOREM 2. - {\it Let $X$ be a topological space and let $I, J:X\to {\bf R}$ be two sequentially
lower semicontinuous functions. Assume that $J$ is sequentially inf-compact and that, for some $c>0$, one has
$$\inf_{x\in J^{-1}(]c,+\infty[)}{{I(x)}\over {J(x)}}=-\infty\ .\eqno{(1)}$$
Then, there exists a divergent sequence $\{\lambda_n^*\}$ in $]0,+\infty[$ with the following property: for every
$n\in {\bf N}$ and for every increasing 
and lower semicontinuous
function $\varphi:J(X)\to {\bf R}$ such that $I+\mu\varphi\circ J$ is sequentially inf-compact for all $\mu>0$,
there exists $\mu^*>0$ such that the function $I+\lambda_n^* J+\mu^*\varphi\circ J$ has at least two global
minima in $X$.}\par
\smallskip
In turn, to prove Theorem 2, we need the two following results that we established in [6] and [7] respectively:\par
\medskip
THEOREM A. - {\it Let $X$ be a topological space and let $\Phi, \Psi:X\to {\bf R}$ be
two functions such that, for every $\lambda>0$, the function  $\Phi+\lambda\Psi$ is sequentially lower semicontinuos and
sequentially inf-compact, and has a unique global minimum in $X$. Assume also that
$\Phi$ has no global minimum.\par
Then, for every $r\in ]\inf_{X}\Psi,\sup_X\Psi[$, there exists $\hat\lambda_r>0$ such that the unique global
minimum in $X$ of the function $\Phi+\hat\lambda_r\Psi$ lies in $\Psi^{-1}(r)$\ .}\par
\medskip
THEOREM B. - {\it 
 Let $S$ be a topological space
 and let $P, Q:S\to {\bf R}$ be two functions satisfying
the following conditions:\par
\noindent
$(a)$\hskip 5pt for each $\lambda>0$, the
function $P+\lambda Q$ is sequentially lower semicontinuous and sequentially inf-compact;\par
\noindent
$(b)$ there exist $\rho\in ]\inf_X Q,\sup_X Q[$ and
$v_1, v_2\in X$ such that
$$Q(v_1)<\rho<Q(v_2)\eqno{(2)}$$
and
$${{P(v_1)-\inf_{Q^{-1}(]-\infty,\rho])}P}\over {\rho-Q(v_1)}}<
{{P(v_2)-\inf_{Q^{-1}(]-\infty,\rho])}P}\over {\rho-Q(v_2)}}\ .\eqno{(3)}$$
Under such hypotheses, there exists $\lambda^*>0$ such that the function
$P+\lambda^*Q$ has at least two global minima.}\par
\medskip
{\it Proof of Theorem 2.}  Fix
 $\rho_0>\inf_XJ$, $x_0\in J^{-1}(]-\infty,\rho_0[)$ 
 and $\lambda$ satisfying
$$\lambda>{{I(x_0)-\inf_{J^{-1}(]-\infty,\rho_0])}I}\over {\rho_0-J(x_0)}}\ .$$
Hence, one has
$$I(x_0)+\lambda J(x_0)<\lambda\rho_0+\inf_{J^{-1}(]-\infty,\rho_0])}I
\ .\eqno{(4)}$$
Since $J^{-1}(]-\infty,\rho_0])$ is sequentially compact, by sequential lower
semicontinuity, there is
$\hat x\in J^{-1}(]-\infty,\rho_0])$ such that
$$I(\hat x)+\lambda J(\hat x)=
\inf_{x\in J^{-1}(]-\infty,\rho_0])}(I(x)+\lambda J(x))\ .\eqno{(5)}$$
We claim that 
$$J(\hat x)<\rho_0\ .\eqno{(6)}$$ 
Arguing by contradiction,
assume that $J(\hat x)=\rho_0$. Then, in view of $(4)$, we would have
$$I(x_0)+\lambda J(x_0)<I(\hat x)+\lambda J(\hat x)$$
against $(5)$. By $(1)$, there is a sequence $\{x_n\}$ in
$J^{-1}(]c,+\infty[)$ such that 
$$\lim_{n\to \infty}{{I(x_n)}\over {J(x_n)}}=-\infty\ .$$
Now, set
$$\gamma=\min\left \{ 0,\inf_{x\in J^{-1}(]-\infty,\rho_0])}(I(x)+\lambda J(x))\right \}$$
and fix $\hat n\in {\bf N}$ so that
$${{I(x_{\hat n})}\over {J(x_{\hat n})}}<-\lambda + {{\gamma}\over {c}}\ .$$
We then have
$$I(x_{\hat n})+\lambda J(x_{\hat n}) < {{\gamma}\over {c}}J(x_{\hat n})\leq \gamma\leq
\inf_{x\in J^{-1}(]-\infty,\rho_0])}(I(x)+\lambda J(x))\ .\eqno{(7)}$$
In particular, this implies that
 $$J( x_{\hat n})>\rho_0\ .\eqno{(8)}$$
Put
$$\rho_{\lambda}^*=J(x_{\hat n})\ .$$
At this point,  we realize that it is possible
to apply Theorem B taking 
$$S=J^{-1}(]-\infty,\rho_{\lambda}^*])\ ,$$ 
$$P=I_{|S}+\lambda J_{|S}\ ,$$ 
$$Q=J_{|S}\ .$$
 Indeed, $(a)$ is satisfied since $S$ is sequentially compact. To satisfy $(b)$, take
$$\rho=\rho_0\ ,$$
$$v_1=\hat x\ ,$$
$$v_2=x_{\hat n}\ .$$
So, with these choices, $(2)$ follows from $(6)$ and $(8)$, while $(3)$ follows from $(5)$ and $(7)$.
Consequently, Theorem B ensures the existence of $\delta_{\lambda}>0$ such that
 the restriction
of the function $I+(\lambda+\delta_{\lambda})J$ to
$J^{-1}(]-\infty,\rho_{\lambda}^*])$ has at least two global minima, say $w_1, w_2$.
Now, fix an increasing and lower semicontinuous function $\varphi:J(X)\to {\bf R}$ such that
$I+\mu\varphi\circ J$ is sequentially inf-compact for all $\mu>0$. We claim that, for some $\mu>0$,
the function $I+(\lambda+\delta_{\lambda})J+\mu\varphi\circ J$ has at least two global minima in $X$.
Arguing by contradiction, assume that, for each $\mu>0$, there exists a unique global minimum in $X$ for the
function $I+(\lambda+\delta_{\lambda})J+\mu\varphi\circ J$ (which is clearly sequentially lower semicontinuous and
sequentially inf-compact). Now, after observing that, by $(1)$, the function $I+(\lambda+\delta_{\lambda})J$ is unbounded below,
we can apply Theorem A taking
$$\Phi=I+(\lambda+\delta_{\lambda})J$$
and
$$\Psi=\varphi\circ J\ .$$ 
Observe that the function $\varphi\circ J$ is unbounded above. Indeed, if not, the sequential inf-compactness of $\varphi\circ J$ jointly with
the sequential lower semicontinuity of $I$ would imply the existence of a global minimum for $I$, contrary to $(1)$.
Moreover, since $J(x_0)<\rho_{\lambda}^*$, we have
$$\inf_X\varphi\circ J\leq \varphi(J(x_0))<\varphi(\rho_{\lambda}^*)\ .$$
Then, Theorem A ensures the existence of $\hat\mu>0$ such that the unique global minimum in $X$
of the function $I+(\lambda+\delta_{\lambda})J+\hat\mu\varphi\circ J$, say $\hat w$, lies in $(\varphi\circ J)^{-1}(\varphi(\rho_{\lambda}^*))$.
Since $\varphi$ is increasing, we have
$$J^{-1}(]-\infty,\rho_{\lambda}^*])=(\varphi\circ J)^{-1}(]-\infty,\varphi(\rho_{\lambda}^*)])$$
and hence, for
$i=1,2$, we have
$$\inf_{x\in X}(I(x)+(\lambda+\delta_{\lambda})J(x)+\hat\mu\varphi(J(x)))\leq 
I(w_i)+(\lambda+\delta_{\lambda})J(w_i)+\hat\mu\varphi(J(w_i))$$
$$\leq I(\hat w)+(\lambda+\delta_{\lambda})J(\hat w)+\hat\mu\varphi(J(\hat w))=\inf_{x\in X}(I(x)+(\lambda+\delta_{\lambda})J(x)+\hat\mu\varphi(J(x)))\ .$$
That is to say, $w_1$ and $w_2$ would be two global minima in $X$ of the function $I+(\lambda+\delta_{\lambda})J+\hat\mu\varphi\circ J$,
a contradiction. Therefore, it remains proved that there exists $\mu^*>0$ such that the function $I+(\lambda+\delta_{\lambda})J+\mu^*\varphi\circ J$
has at least two global minima in $X$. Finally, observe that the set
$$A:=\left \{\lambda+\delta_{\lambda} : \lambda>{{I(x_0)-\inf_{J^{-1}(]-\infty,\rho_0])}I}\over {\rho_0-J(x_0)}}\right \}$$
is unbounded above. So, for what we have seen above, any divergent sequence $\{\lambda_n^*\}$ in $A$ satisfies the thesis.\hfill
$\bigtriangleup$\par
\medskip
{\it Proof of Theorem 1}.  Fix $q\in \left ] 0, {{m+2}\over {m-2}}\right [$ and $f\in {\cal A}_q$.
We are going to apply Theorem 2 taking $X=H^1_0(\Omega)$, endowed with the weak topology, and
$I, J:H^1_0(\Omega)\to {\bf R}$ defined by 
$$I(u)=-\int_{\Omega}F(u(x))dx\ ,$$
$$J(u)={{1}\over {2}}\|u\|^2\ ,$$
where
$$\|u\|^2=\int_{\Omega}|\nabla u(x)|^2dx\ .$$
Clearly,  $J$ is weakly inf-compact and $I$ (since $f$ has a subcritical growth) is sequentially weakly continuous. Now, fix a measurable set
$C\subset\Omega$, of positive measure, and a function $w\in H^1_0(\Omega)$ such that $w(x)=1$ for all $x\in C$. Since $f\in {\cal A}_q$, there
exist a sequence $\{\xi_n\}$ in ${\bf R}$, with $\lim_{n\to \infty}|\xi_n|=+\infty$, and a constant $\alpha>0$ such that
$$-\alpha(\xi^2+1)\leq F(\xi)$$
for all $\xi\in {\bf R}$ and
$$\lim_{n\to +\infty}{{F(\xi_n)}\over {\xi_n^2}}=+\infty\ .$$
Thus, we have
$${{\int_{\Omega}F(\xi_nw(x))dx}\over {\int_{\Omega}|\nabla \xi_nw(x)|^2dx}}
={{\hbox {\rm meas}(C)F(\xi_n)+\int_{\Omega\setminus C}F(\xi_nw(x)dx}
\over {\xi_n^2\int_{\Omega}|\nabla w(x)|^2dx}}$$
$$\geq {{\hbox {\rm meas}(C)F(\xi_n)}\over {\xi_n^2\int_{\Omega}|\nabla w(x)|^2dx}}-\alpha{{\int_{\Omega}|w(x)|^2dx+{{\hbox {\rm meas}(\Omega)}
\over {\xi_n^2}}}\over {\int_{\Omega}|\nabla w(x)|^2dx}}$$
and so
$$\liminf_{\|u\|\to +\infty}{{I(u)}\over {J(u)}}=-\infty\ .$$
Therefore, the assumptions of Theorem 2 are satisfied. Let $\{\lambda_n^*\}$ be a divergent sequence with the property expressed in Theorem 2.
Fix $n\in {\bf N}$ and a continuous and
non-decreasing function $k:[0,+\infty[\to [0,+\infty[$, with $\lim_ {t\to +\infty}{{K(t)}\over {t^{q+1\over 2}}}=+\infty$ and int$(k^{-1}(0))=\emptyset$. Let $\varphi:[0,+\infty[\to [0,+\infty[$ be defined by
$$\varphi(t)={{1}\over {2}}K(2t)$$
for all $t\geq 0$. Clearly, the function $\varphi$ is increasing (and continuous).  Moreover, due to the Sobolev imbedding, there is a constant
$\beta>0$ such that
$$I(u)\geq -\beta(1+\|u\|^{q+1})$$
for all $u\in X$ and so, for each $\mu>0$, we have
$$I(u)+\mu\varphi(J(u))\geq -\beta(1+\|u\|^{q+1})+{{\mu}\over {2}}K(\|u\|^2)=
\|u\|^{q+1}\left ( -\beta\left ( 1+{{1}\over {\|u\|^{q+1}}}\right ) +{{\mu}\over {2}}{{K(\|u\|^2)}\over {\|u\|^{q+1}}}\right )\eqno{(9)}$$
for all $u\in X$. Since
$$\lim_{\|u\|\to +\infty}{{K(\|u\|^2)}\over {\|u\|^{q+1}}}=+\infty\ ,$$
from $(9)$ we infer that the functional $I+\mu\varphi\circ J$ is sequentially weakly inf-compact. As a consequence, there exists $\mu^*>0$ such that
the functional $I+\lambda_n^*J+\mu^*\varphi\circ J$ has at least two global minima in $X$ which, therefore, are
weak solutions of the problem we are dealing with. Now,
  observe that the function $t\to t(\lambda_n^*+\mu^*k(t^2))$  is
increasing in $[0,+\infty[$ and its range is $[0,+\infty[$. Denote by $\psi$ its inverse.
Let $T:X\to X$ be the operator defined by
$$T(v)=\cases {{{\psi(\|v\|)}\over {\|v\|}}v & if $v\neq 0$\cr & \cr
0 & if $v=0$\ ,\cr}$$
 Since $\psi$ is continuous
and $\psi(0)=0$, the operator $T$ is continuous in $X$. For each $u\in X\setminus \{0\}$,
we have
$$T((\lambda_n^*+\mu^*k(\|u\|^2))u)={{\psi((\lambda_n^*+\mu^*k(\|u\|^2))\|u\|)}\over {(\lambda_n^*+\mu^*k(\|u\|^2))\|u\|}}
(\lambda_n^*+\mu^*k(\|u\|^2))u$$
$$={{\|u\|}\over {(\lambda_n^*+\mu^*k(\|u\|^2))\|u\|}}(\lambda_n^*+\mu^*k(\|u\|^2))u=u\ .$$
In other words, $T$ is a continuous inverse of the derivative of the functional $\lambda_n^*J+\mu^*\varphi\circ J$.
Then, since the derivative of $I$ is compact, the functional $I+\lambda_n^*J+\mu^*\varphi\circ J$
satisfies the Palais-Smale condition ([10], Example 38.25) and hence the existence of a third critical
point of the same functional is assured by Corollary 1 of [3]. The proof is complete.\hfill $\bigtriangleup$\par
\medskip
We conclude by formulating two open problems.\par
\medskip
PROBLEM 1. - In Theorem 1, can the role of the sequence $\{a_n\}$ be assumed by a suitable unbounded interval ?\par
\medskip
PROBLEM 2. - Does Theorem 1 hold for $q={{m+2}\over {m-2}}$ ?
\bigskip
\bigskip
{\bf Acknowledgement.} The author has been supported by the Gruppo Nazionale per l'Analisi Matematica,
la Probabilit\`a e le loro Applicazioni (GNAMPA) of the Istituto Nazionale di Alta Matematica (INdAM).\par
\vfill\eject
\bigskip
\bigskip
\centerline {\bf References}\par
\bigskip
\bigskip
\bigskip
\bigskip
\noindent
[1]\hskip 5pt G. M. FIGUEIREDO and R. G. NASCIMENTO, {\it Existence of a nodal solution with minimal energy for a Kirchhoff equation},
Math. Nachr., {\bf 288} (2015), 48-60.\par
\smallskip
\noindent
[2]\hskip 5pt K. PERERA and Z. T. ZHANG, {\it Nontrivial solutions of
Kirchhoff-type problems via the Yang index}, J. Differential Equations,
{\bf 221} (2006), 246-255.\par
\smallskip
\noindent
[3]\hskip 5pt P. PUCCI and J. SERRIN, {\it A mountain pass theorem},
J. Differential Equations, {\bf 60} (1985), 142-149.\par
\smallskip
\noindent
[4]\hskip 5pt B. RICCERI, {\it On an elliptic Kirchhoff-type problem
depending on two parameters}, J. Gobal Optim., {\bf 46}
(2010), 543-549.\par
\smallskip
\noindent
[5]\hskip 5pt B. RICCERI, {\it Energy functionals of Kirchhoff-type problems having multiple global minima}, Nonlinear Anal., {\bf 115} (2015),  
130-136.\par
\smallskip
\noindent
[6]\hskip 5pt B. RICCERI, {\it Well-posedness of constrained minimization
problems via saddle-points}, J. Global Optim., {\bf 40} (2008),
389-397.\par
\smallskip
\noindent
[7]\hskip 5pt B. RICCERI, {\it Multiplicity of global minima for
parametrized functions}, Rend. Lincei Mat. Appl., {\bf 21} (2010),
47-57.\par
\smallskip
\noindent
[8]\hskip 5pt J. SUN and T. F. WU, {\it Existence and multiplicity of solutions for an indefinite Kirchhoff-type equation in bounded domains}, 
Proc. Roy. Soc. Edinburgh Sect. A,  {\bf 146} (2016), 435-448. \par
\smallskip
\noindent
[9]\hskip 5pt X. H. TANG and B. CHENG, {\it Ground state sign-changing solutions for Kirchhoff type problems in bounded domains}, J. Differential Equations, {\bf  261} (2016), 2384-2402.\par
\smallskip
\noindent
[10] E. ZEIDLER, {\it Nonlinear functional analysis and its applications}, vol. III, Springer-Verlag, 1985.\par
\smallskip
\noindent
[11]\hskip 5pt Q. G. ZHANG, H. R. SUN and J. J. NIETO, {\it J. Positive solution for a superlinear Kirchhoff type problem with a parameter}, Nonlinear Anal.,
 {\bf 95} (2014), 333-338. 

\bigskip
\bigskip
\bigskip
\bigskip
Department of Mathematics\par
University of Catania\par
Viale A. Doria 6\par
95125 Catania\par
Italy\par
{\it e-mail address:} ricceri@dmi.unict.it

\bye